\documentclass[11pt,twoside,english]{myamsart}
\advance\oddsidemargin by -1.0cm
\advance\evensidemargin by -1.0cm
\advance\topmargin by -1.0cm 
\usepackage{amssymb}
\usepackage{babel}
\usepackage{amstext}
\usepackage{amscd}   
\usepackage{epsfig}  
\usepackage{rotating}
\let\mathg\mathfrak
\theoremstyle{plain}
\newtheorem{cor}{Corollary}[section]

\newtheorem{thm}{Theorem}[section]            
\newtheorem{prop}{Proposition}[section]
\theoremstyle{definition}
\newtheorem{exa}{Example}[section]
\newtheorem{NB}{Remark}[section]
\newtheorem{dfn}{Definition}[section]

\renewcommand{\labelenumi}{$(\theenumi)$}
\newcommand{\bdm}{\begin{displaymath}}
\newcommand{\edm}{\end{displaymath}}
\newcommand{\be}{\begin{equation}}
\newcommand{\ee}{\end{equation}}
\newcommand{\ba}[1]{\begin{array}{#1}}
\newcommand{\ea}{\end{array}}

\newcommand{\btab}{\begin{tabular}}
\newcommand{\etab}{\end{tabular}}

\newcommand{\R}{\ensuremath{\mathbb{R}}}

\newcommand{\T}{\ensuremath{\mathrm{T}}}

\newcommand{\G}{\ensuremath{\mathrm{G}}}
\newcommand{\J}{\ensuremath{\mathrm{J}}}

\newcommand{\Ric}{\ensuremath{\mathrm{Ric}}}
\newcommand{\Scal}{\ensuremath{\mathrm{Scal}}}

\newcommand{\U}{\ensuremath{\mathrm{U}}}

\newcommand{\so}{\ensuremath{\mathg{so}}}
\newcommand{\SO}{\ensuremath{\mathrm{SO}}}
\newcommand{\Spin}{\ensuremath{\mathrm{Spin}}}

\newcommand{\g}{\ensuremath{\mathfrak{g}}}

\newcommand{\m}{\ensuremath{\mathfrak{m}}}

\begin{document}
\def\haken{\mathbin{\hbox to 6pt{%
                 \vrule height0.4pt width5pt depth0pt
                 \kern-.4pt
                 \vrule height6pt width0.4pt depth0pt\hss}}}
    \let \hook\intprod
\setcounter{equation}{0}
%
%
\thispagestyle{empty}
%
\date{\today}
\title[Geometric structures of vectorial type]
{Geometric structures of vectorial type}
%
%
%
\author{I. Agricola}
\author{T. Friedrich}
\address{\hspace{-5mm} 
{\normalfont\ttfamily agricola@mathematik.hu-berlin.de\newline
\normalfont\ttfamily friedric@mathematik.hu-berlin.de}\newline
Institut f\"ur Mathematik \newline
Humboldt-Universit\"at zu Berlin\newline
Sitz: WBC Adlershof\newline
D-10099 Berlin, Germany}
%
\thanks{Supported by the SFB 647 "Raum, Zeit,Materie" and the SPP 1154 
``Globale Differentialgeometrie'' of the DFG as well as the
Volkswagen Foundation}
\subjclass[2000]{Primary 53 C 25; Secondary 81 T 30}
\begin{abstract}
We study geometric structures of $\mathcal{W}_4$-type in the sense of A. Gray
on a Riemannian manifold. If the structure group $\mathrm{G} \subset \SO(n)$
preserves a spinor or a non-degenerate differential form,  its intrinsic
torsion $\Gamma$ is a closed $1$-form (Proposition \ref{dGamma} and Theorem
\ref{Fixspinor}). Using a $\mathrm{G}$-invariant spinor we prove a splitting
theorem (Proposition \ref{splitting}). The latter result generalizes and 
unifies a recent result obtained in \cite{Ivanov&Co}, 
where this splitting has been
proved in dimensions $n=7,8$ only. Finally we investigate geometric
structures of vectorial type and admitting a characteristic connection
$\nabla^{\mathrm{c}}$.
An interesting class of geometric structures 
generalizing Hopf structures are those with a
$\nabla^{\mathrm{c}}$-parallel intrinsic torsion $\Gamma$. 
In this case, $\Gamma$ induces a Killing vector
field (Proposition \ref{Killing}) and  
for some special structure groups it is even parallel.
\end{abstract}
\maketitle
\pagestyle{headings}
%
%
%
\section{Adapted connections of a geometric structure
of vectorial type}\noindent

\noindent
Fix a subgroup $\mathrm{G} \subset \SO(n)$ of the special
orthogonal group and decompose the Lie algebra $\so(n) = \g \oplus \m$ 
into the Lie algebra $\g$ of $\mathrm{G}$ and its orthogonal complement $\m$. 
The different geometric types of  
$\mathrm{G}$-structures on a Riemannian manifold correspond to the 
irreducible $\mathrm{G}$-components of the representation $\R^n \otimes \m$.
Indeed, consider an oriented Riemannian manifold 
$(M^n, g)$ and denote 
its Riemannian frame bundle by 
$\mathcal{F}(M^n)$. It is a principal $\SO(n)$-bundle 
over $M^n$. A $\mathrm{G}$-structure is a reduction 
$\mathcal{R} \subset \mathcal{F}(M^{n})$ of the frame bundle to the subgroup 
$\mathrm{G}$. The Levi-Civita connection is a $1$-form $Z$ on 
$\mathcal{F}(M^{n})$ with values in the Lie algebra $\mathg{so}(n)$.
We restrict the Levi-Civita connection to $\mathcal{R}$ and decompose 
it with respect to the decomposition of the Lie algebra $\mathg{so}(n)$,
\bdm
Z\big|_{T(\mathcal{R})} \ := \ Z^* \, \oplus \ \Gamma \, .
\edm
Then, $Z^*$ is a connection in the principal $\mathrm{G}$-bundle $\mathcal{R}$
and $\Gamma$ is a $1$-form on $M^{n}$ 
with values in the associated bundle $\mathcal{R} \times_{\mathrm{G}} 
\mathg{m}$. This $1$-form, or more precisely the $\G$-components
of the element $\Gamma \in \R^n \otimes \mathg{m}$,
characterizes the different types of non-integrable
$\G$-structures (see \cite{Fri2}).
The $1$-form $\Gamma$ is called the {\it intrinsic torsion}
of the $\mathrm{G}$-structure. There is a second notion for 
$\mathrm{G}$-structures, namely that of {\it characteristic connection} 
and its {\it characteristic torsion}. By definition,
a characteristic connection is a
$\mathrm{G}$-connection $\nabla^{\mathrm{c}}$ with totally skew symmetric 
torsion tensor. Typically, not every type of $\mathrm{G}$-structure admits a 
characteristic connection. In order to formulate the condition, we embed 
the space of all $3$-forms into $\R^n \otimes \m$ using
the morphism
\bdm
\Theta \, : \, \Lambda^3(\R^n) \longrightarrow \R^n \otimes \m \, , \quad
\Theta(\T) \ := \ \sum_{i=1}^n e_i \otimes \mathrm{pr}_{\m}(e_i \haken
\T) \, .
\edm
A $\mathrm{G}$-structure admits a characteristic connection 
$\nabla^{\mathrm{c}}$ if
and only if the intrinsic torsion $\Gamma$ belongs to the image of the 
$\Theta$. In this case, the characteristic torsion is the pre-image of
the intrinsic torsion (see \cite{Fri2})
\bdm
2 \, \Gamma \ = \ - \,  \Theta(\T^{\mathrm{c}}) \,  \quad \mathrm{and} \quad
\ \nabla^{\mathrm{c}}_X Y \ = \ \nabla^g_X Y \, + \, \frac{1}{2}
\T^{\mathrm{c}}(X\, ,\, Y \, , \, - ) \, .
\edm
For different geometric structures, the characteristic torsion form 
has been computed 
explicitly in terms of the underlying geometric data. Formulas of that type
are known for almost hermitian structures, almost
contact metric structures and $\mathrm{G}_2$- and $\Spin(7)$-structures 
in dimensions seven and eight. 
In case of a Riemannian naturally reductive
space $M^n = \mathrm{G}_1/\mathrm{G}$, we obtain a $\mathrm{G}$-reduction
$\mathcal{R} := \mathrm{G}_1\subset \mathcal{F}(M^{n})$ of the frame bundle.
Then the characteristic connection of the $\mathrm{G}$-structure
coincides with the {\it canonical} connection of the reductive space.
In this sense, we can understand the characteristic connection of a Riemannian 
$\mathrm{G}$-structure as a generalization of the canonical connection
of a Riemannian naturally reductive space. The canonical 
connection of a naturally reductive space has parallel torsion form and
parallel curvature tensor,
$\nabla^{\mathrm{c}} \T^{\mathrm{c}} =  0 =
\nabla^{\mathrm{c}} \mathrm{R}^{\mathrm{c}} $. 
For arbitrary $\mathrm{G}$-structures and their characteristic connections,
these properties do not hold anymore. Corresponding examples are discussed
in \cite{FriedrichIvanov}.
The space $\R^n \otimes \m$ contains $\R^n$ in a natural way,
\bdm
\Theta_1 \, : \, \R^n \longrightarrow \R^n \otimes \m \, , \quad
\Theta_1(\Gamma) \ = \ \sum_{i = 1}^n e_i \otimes
\mathrm{pr}_{\m}(e_i \wedge \Gamma) \ .
\edm
The class of geometric structures we will study in this paper is 
the following one.

\begin{dfn}
Let $M^n$ be an oriented Riemannian manifold and denote by $\mathcal{F}(M^n)$ 
its
frame bundle. A geometric structure $\mathcal{R} \subset \mathcal{F}(M^n)$ 
is called
of {\it vectorial type} if its intrinsic torsion belongs to
$\Gamma \in \R^n \subset \R^n \otimes \m$. 
\end{dfn}

\begin{NB}
These geometric structures are  usually called 
$\mathcal{W}_4$-structures. They occur in the description
of almost hermitian manifolds, of
$\G_2$-structures in dimension seven, of $\Spin(7)$-structures
in dimension eight and $\Spin(9)$-structures
in dimension sixteen (see \cite{Fri0}). 
\end{NB}

\noindent
We will identify vectors field on the Riemannian manifold
$(M^n,g)$ with $1$-forms. Denoting 
the vector field corresponding to the intrinsic
torsion by $\Gamma$, too, we obtain the following formula
for the intrinsic torsion defined  by a 
vector field $\Gamma$,
\bdm
\Gamma(X) \ = \ \mathrm{pr}_{\mathg{m}}(X \wedge \Gamma) \, .
\edm
Of course, a geometric structure of vectorial type does not have
to admit a characteristic connection. It depends on the decomposition of
the $\G$-representation $\Lambda^3(\R^n)$. For example, consider
the subgroup $\SO(3) \subset \SO(5)$ defined by the $5$-dimensional,
real representation of $\SO(3)$. Then $\Lambda^3(\R^5)$ splits into
a $3$-dimensional and a $7$-dimensional $\SO(3)$-representation,
i.\,e., $\Theta(\Lambda^3(\R^5)) \subset \R^5 \otimes \mathg{m}$ and
$\R^5 \subset \R^5 \otimes \mathg{m}$ are complementary subspaces. 
A similar situation occurs for the subgroups $\Spin(9) \subset \SO(16)$
and for $\mathrm{F}_4 \subset \SO(26)$ (see \cite{Fri2}). On the other side,
many interesting geometric structures of vectorial type admit
characteristic connection. This situation occurs
for example for the subgroups $\U(n) \subset \SO(2n)$, $\G_2 \subset
\SO(7)$ and $\Spin(7) \subset \SO(8)$. The corresponding characteristic
torsion has been computed explicitely for these cases in 
\cite{FriedrichIvanov}.\\

\noindent
The first observation is a link to E.~Cartan (see \cite{AgFr1}, 
\cite{Cartan1}), who classified
the types of metric connections. There are two special classes. The first
class are metric connections of vectorial type, the second class are metric
connections with a totally skew-symmetric torsion. There also exists a third
class, but does not have direct geometric interpretation.
The geodesic flow of metric connections
of vectorial type has been investigated in \cite{Cartan1} and \cite{AgTh}. 
On the other side, the geodesic flow of metric connections with totally
skew-symmetric torsion coincides with the Riemannian geodesic flow.

\begin{prop}
If a $\G$-structure is of vectorial type, then there exists a unique
metric connection $\nabla^{\mathrm{vec}}$ 
of vectorial type in the sense of Cartan and preserving
the $\G$-structure. The formula is
\bdm
\nabla^{\mathrm{vec}}_XY\ =\ \nabla^g_XY\,-\,g(X \, , \,Y)\cdot\Gamma\, + \, 
g(Y \, ,\,\Gamma)\cdot X \, .
\edm
Conversely, if a $\G$-structure $\mathcal{R}$ admits a connection of 
vectorial type in the sense of Cartan, then $\mathcal{R}$ is of vectorial type
in our sense.
\end{prop}
\begin{proof}
The Levi-Civita connections splits into 
\bdm
Z(X) \ = \ Z^*(X) \, + \,
\mathrm{pr}_{\mathg{m}}(X \wedge \Gamma) = Z^*(X) \, - \, \mathrm{pr}_{\g}(X
\wedge \Gamma)
\, + \, X \wedge \Gamma \, .
\edm
The formula $\beta(X):= \mathrm{pr}_{\g}(X \wedge \Gamma)$ defines a $1$-form
with values
in the Lie algebra $\g$. Therefore, the connection
\bdm
Z^{\mathrm{vec}}(X) \ := \ Z(X) \, - \, X \wedge \Gamma \ = \ 
Z^*(X) \, - \,\mathrm{pr}_{\g}(X \wedge \Gamma) 
\edm 
is a $\G$-connection. It is of vectorial type in the sense of Cartan. 
Suppose vice versa that there exists a $\G$-connection $Z^{**}$ of
vectorial type. We compare it with the Levi-Civita connection and obtain
the relation
\bdm
Z^{**}(X)  \ = \ Z(X) \, + \, X \wedge \Gamma \, .
\edm
Moreover, the definition of the $1$-form $\Gamma$ as well as the 
$\mbox{G}$-connection $Z^*$ yields the equation
\bdm
Z(X) \ = \ Z^*(X) \, + \, \Gamma(X) \, .
\edm
Finally, since $Z^{**}$ preserves the $\mbox{G}$-structure, there
exists a $1$-form $\beta$ with values in the Lie algebra $\mathg{g}$ such that
\bdm
Z^{**}(X)  \ = \ Z^*(X)  \, + \, \beta(X) \, .
\edm
Combining these three formulas we obtain, for any vector $X$, the equation
\bdm
\beta(X) \ = \ \Gamma(X) \, + \, X \wedge \Gamma \, .
\edm
We now take the projection onto $\mathg{m}$. Since $\beta(X)$ belongs
to the Lie algebra $\mathg{g}$, we conclude that $\Gamma$ should be in the 
image of $\R^n \subset  \R^n \otimes \m$.
\end{proof}
\noindent
Let $\Gamma$ be a vector field on a Riemannian manifold $(M^n,g)$.
Then we define a metric connection $\nabla^{\mathrm{vec}}$ as before.
Its holonomy group is a subgroup of $\SO(n)$ and its holonomy
bundle is a reduction of the frame bundle $\mathcal{F}(M^n)$. Obviously,
the corresponding structure is of vectorial type. Therefore, any vector
field can occur. However, if the structure group $\G$ is fixed, then
we obtain restrictions for the possible vector field $\Gamma$. In the next
section we will explain the corresponding results.
\vspace{2mm}

\noindent
A geometric structure of vectorial type induces a triple $(M^n, g, \Gamma)$
consisting of a Riemannian manifold and a vector field. A similar situation
occurs in Weyl geometry (see \cite{Calderbank}, \cite{Gauduchon}).
A Weyl structure is a pair consisting of a conformal
class of metric and a torsion free connection preserving 
the conformal structure.
Choosing a metric $g$ in the conformal class, the connection defines 
a vector field and the corresponding covariant derivative on vectors
is defined by the formula
\bdm
\nabla^{\mathrm{w}}_X Y\ =\ \nabla^g_X Y\,+\,g(X\, ,\,\Gamma)\cdot Y\,+\, 
g(Y \, , \, \Gamma) \cdot X \, - \, g(X \, , \, Y) \cdot \Gamma \, .
\edm
Weyl geometry deals with the geometric properties of these connections. 
The two connections $\nabla^{\mathrm{vec}}$ and $\nabla^{\mathrm{w}}$ 
are different. The Weyl connection does not preserve 
any Riemannian geometric structure, moreover, it is torsion free.
However, the curvature tensors and the Ricci tensors are closely
related,
5
\bdm
\mathcal{R}^{\mathrm{vec}}(X , Y) Z \ = \ \mathcal{R}^{\mathrm{w}}(X,Y)Z \,
- \, d \Gamma(X , Y) \cdot Z \, , \quad
\Ric^{\mathrm{vec}} \ = \ \Ric^{\mathrm{w}} \, + \ d\Gamma \, .
\edm 
In particular, the symmetric parts of the Ricci tensors coincide.
If one can prove that 
certain geometric structures induce  Weyl-Einstein structures,
one can apply several results known in Weyl geometry. 
Examples of this approach can be found in Theorem \ref{Fixspinor} and  
Proposition \ref{splitting}. Otherwise, the topics are quite different.     
\vspace{2mm}

\noindent
Of course, a conformal change of $\G$-structures is again possible.
Let us discuss it. The total space
$\mathcal{R} \subset \mathcal{F}(M^n,g)$ of a $\G$-structure consists
of $n$-tuples $(e_1, e_2 , \ldots , e_n)$ of orthonormal vectors
tangent to $M^n$. Let $g^* := e^{2f} g$ be a conformal change
of the metric. Then we define a new $\G$-structure
$\mathcal{R}^* \subset \mathcal{F}(M^n, g^*)$ by
\bdm
\mathcal{R}^* \ = \ \Big\{(e^{-f} \cdot e_1\ , \ e^{-f} \cdot e_2 \
\ldots \ , \ e^{-f} \cdot e_n) \ : \ (e_1 \ , \ e_2 \ , \ \ldots \ , \ e_n) 
\in  \mathcal{R} \Big\} \, .
\edm
The intrinsic torsion changes by the
element $df \in \R^n \subset \R^n \otimes \mathg{m}$,
$\Gamma^* =  \Gamma  +  df$.
In particular, the conformal change of a geometric structure of vectorial
type is again of vectorial type. Moreover, the differentials
$d \Gamma = d \Gamma^*$ coincide. On the other side, starting
with an arbitrary geometric structure on a compact manifold, the equation
\bdm
0 \ = \ \delta^{g^*}( \Gamma^*) \ = \ \delta^g(\Gamma) \, + \, \Delta(f) \, 
+ \, (n - 2) \cdot \big( (df , \Gamma) \, + \, ||df||^2 \big)
\edm
has a unique solution $f = - \, \Delta^{-1} ( \delta^g(\Gamma))$. 
Consequently, an arbitrary geometric
structure of vectorial type on a compact manifold
admits a conformal change such that 
the new $1$-form is coclosed (see \cite{Gauduchon}, \cite{Calderbank}). 
In principle,
one can reduce the investigation of geometric structures of vectorial type
on compact manifolds to those structures with a coclosed form,  
$\delta^g(\Gamma) = 0$.
%
\section{Parallel forms and spinors}\noindent
%
\noindent
Let $\Omega^k \in \Lambda^k(\R^n)$ be a
$\G$-invariant $k$-form. It defines a $k$-form on any Riemannian manifold
with a fixed $\G$-structure which is parallel with respect
to any $\G$-connection. The Lie algebra $\so(n) = \Lambda^2(\R^n)$
acts on the vectors space $\Lambda^k(\R^n)$ via the formula
\bdm
\rho_*(\omega^2) \, \Omega^k \ = \ \sum_{i = 1}^n (e_i \haken \omega^2) \wedge
(e_i \haken \Omega^k) \, , \quad \omega^2  \in \ \so(n) \, .
\edm
Consequently, we can compute the Riemannian covariant derivative of
$\Omega^k$,
\begin{eqnarray*}
\nabla^g_X \Omega^k &= & \rho_* \big(\mathrm{pr}_{\mathg{m}}(X \wedge
\Gamma)\big) \, \Omega^k \ = \
\rho_* \big(X \wedge \Gamma\big) \, \Omega^k \ = \
\sum_{i = 1}^n \big( e_i \haken (X \wedge \Gamma)\big)
\wedge 
\big( e_i \haken \Omega^k \big) \\
&=& \Gamma \wedge \big(X \haken \Omega^k \big) \, - \, X \wedge \big(\Gamma
\haken
\Omega^k \big) \, .
\end{eqnarray*}
The differential of $\Omega^k$ as well as its codifferential are given
by
\begin{eqnarray*}
d  \Omega^k &=& \sum_{i=1}^n \, e_i \wedge \nabla^g_{e_i} \, \Omega^k \ = \ -
\, k \cdot (\Gamma \wedge \Omega^k) \, , \\
\delta^{g}  \Omega^k &=& - \, \sum_{i=1}^n \, e_i \haken \nabla^g_{e_i} \, \Omega^k \ = \ (n \, - \, k) \cdot (\Gamma \haken \Omega^k) \ .
\end{eqnarray*}
In particular, for any geometric structure of vectorial type and any
$\G$-invariant form $\Omega^k$ we have
\bdm
\nabla^g_{\Gamma} \, \Omega^k \ = \ 0 \ \quad \mbox{and} \quad \ d \Gamma
\wedge \Omega^k \ = \ 0 \, .
\edm
From these equations we see  that  $\Gamma$ is automatically closed if the 
$k$-form $\Omega^k$
-- treated as a linear map defined on $2$-forms -- has trivial kernel.  
\vspace{2mm}

\begin{prop} \label{dGamma}
Let $\G \subset \SO(n)$ be a subgroup such that
\begin{enumerate}
\item there exists a $\G$-invariant 
differential form $\Omega^k$ of some degree $k$, and
\item the multiplication $\Omega^k \, : \, \Lambda^2(\R^n) \rightarrow
\Lambda^{k+2}(\R^n)$ is injective.
\end{enumerate}

Then, for any $\G$-structure of vectorial type, the $1$-form $\Gamma$ is 
closed, $d \Gamma = 0$.
\end{prop}
\vspace{2mm}

\begin{NB}
The groups $\G_2 \subset \SO(7)$ and $\Spin(7) \subset \SO(8)$ satisfy the
conditions of the Proposition. Consequently, it generalizes results
of Cabrera (see \cite{Cabrera1}, \cite{Cabrera2}). Moreover, there are
other groups satisfying the conditions, namely 
$\U(n) \subset \SO(2n)$ for $n >2$ and $\Spin(9) \subset \SO(16)$. In 
particular, there is an analogon of Cabrera's result for $\Spin(9)$. 
On the other side, there are interesting
$\G$-structures where the group does not satisfy the conditions. The first
example $\SO(3) \subset \SO(5)$ (the irreducible representation)
does not admit any invariant differential form. The subgroup
$\U(2) \subset \SO(4)$ admits an invariant form, but the 
second condition of the Proposition is not satisfied. In these geometries
the condition $d \Gamma = 0$ is an additional requirement on the geometric
structure of vectorial type.
\end{NB}
\vspace{2mm}

\begin{exa} \label{ExaU(2)}
Consider the subgroup $\U(2) \subset \SO(4)$. There are only two types
of $\U(2)$-structures. An almost hermitian manifold
$(M^4,g,\J)$ is of vectorial type
if and only if the almost complex structure is integrable, see 
\cite{AlFrSchoe}. Consequently, starting with an arbitrary 
complex $4$-manifold $(M^4,\J)$, $\it{any}$ hermitian metric $g$ yields
an $\U(2)$-structure of vectorial type and the
vector field $\Gamma$ is defined by the formulas
\bdm
d \Omega \ = \ - \, 2 \, \Gamma \wedge \Omega \ , \quad
\delta^g \Omega \ = \ 2 \, \Gamma \haken \Omega \, .
\edm
Solving this algebraic equation, we obtain
$2 \,\Gamma =  * \, \J(d \Omega)$.
In general, this
$1$-form is not closed. Hermitian manifolds with a closed form $\Gamma$
are called locally conformal
K\"ahler manifolds. In higher dimensions (i.e. for $\G = \U(n)$ and
$n \geq 3$) all hermitian manifolds of vectorial type are automatically
locally conformal K\"ahler. 
\end{exa}
\vspace{2mm}

\begin{exa}
Consider the subgroup $\G = \SO(n-1) \subset \SO(n)$. A $\G$-structure
on $(M^n,g)$ is a vector field $\Omega$ 
of length one. 
The geometric structure is of vectorial type
if and only if there exists a vector field $\Gamma$ such that
\bdm
0 \ = \ \nabla^{\mathrm{vec}}_X \Omega \ = \ 
\nabla^g_X \Omega \, - \, g(X \, , \, \Omega) \, \Gamma \, + \,
g (\Omega \, , \, \Gamma) \, X
\edm
holds. This condition implies that $\Omega$ defines a codimension 
one foliation  on $M^n$,
\bdm
d \Omega \ = \ \Omega \wedge \Gamma \ .
\edm
Moreover,
the second fundamental form of any leave $F^{n-1} \subset M^n$ is given
by the formula
\bdm
\mathrm{II} (X) \ = \ - \, g(\Omega \, , \, \Gamma) \cdot 
X  \, , \quad X \in TF^{n-1} \ .
\edm
Therefore, the leaves are umbilic. Conversely, let $\Omega$ be
a $1$-form defining an umbilic foliation. Let us define the vector field
$\Gamma$ by the formulas
\bdm
\mathrm{II} (X) \ = \ - \, g(\Omega \, , \, \Gamma) \cdot X  \, , \quad \Gamma \ = \ \nabla^g_{\Omega} \Omega \, + \, 
g(\Omega \, , \, \Gamma) \cdot \Omega \ .
\edm
Then the induced $\SO(n-1)$-structure is of vectorial type and
$\Gamma$ is the corresponding vector field. In consequence,
{\it $\SO(n-1)$-structures
of vectorial type coincide with umbilic foliations of codimension one}.
The vector field $\Gamma$ satisfies the condition $\Omega \wedge d \Gamma =
0$, but in general it does not have to be closed.
\end{exa}
\vspace{2mm}

\begin{NB} An almost contact metric structure $(M^{2k + 1}, \xi, \eta,
  \varphi)$ is never of vectorial type. Indeed, the condition 
$\eta \wedge (d\eta)^k \neq 0$ contradicts $ d\eta = \eta \wedge\Gamma$.
\end{NB}
\vspace{2mm}

\noindent
Fix a spin structure of the manifold $(M^n,g)$. The metric connection
$\nabla^{\mathrm{vec}}$ acts on arbitrary spinor fields by
\bdm
\nabla^{\mathrm{vect}}_X \Psi \ = \ \nabla^g_X \Psi \, - \, \frac{1}{2} \cdot
(X \wedge \Gamma) \cdot \Psi \, .
\edm
We now consider the case that the group $\G$ lifts into the spin
group $\Spin(n)$ and admits a $\G$-invariant algebraic spinor 
$\Psi \in \Delta_n$
in the $n$-dimensional spin representation $\Delta_n$. 
We normalize the length of the spinor, $||\Psi|| = 1$.
It defines a spinor field on any Riemannian manifold with a 
$\G$-structure. Moreover, $\Psi$ is parallel with respect to
any $\G$-connection. Using this parallel 
spinor field we can calculate the Riemannian
Ricci tensor $\Ric^g$ completely. Furthermore, we obtain
an algebraic restriction for the $2$-form $d \Gamma$.
\vspace{2mm}

\begin{thm} \label{Fixspinor}
Let $\G \subset \SO(n)$ be a subgroup lifting into the spin group and
suppose that there exists a $\G$-invariant spinor $0 \neq \Psi \in \Delta_n$.
Then the Clifford product
$d  \Gamma \cdot \Psi = 0$ vanishes for any $\G$-structure of vectorial
type. If the dimension $n \geq 5$ is at least five, then 
$\Gamma$ is closed, $d \Gamma = 0$.
The Ricci tensor is given in dimension $n =4$ by 
\begin{eqnarray*}
g(\Ric^g(X)\, , \, Y) \ = \ g(\nabla^g_X\Gamma\, , \, Y)  \, + \, 
 g(\nabla^g_Y\Gamma\, , \, X) \, - \, \delta^g(\Gamma) \cdot g(X \, , \, Y) 
\, + \,   g(A(X \, , \, \Gamma) \, , \, Y)  \, ,
\end{eqnarray*}
and in higher dimensions $n \geq 5$ by
\bdm
\Ric^g(X)\ =\ (n - 2)\,\nabla^g_X\Gamma\,-\,\delta^g(\Gamma) \cdot X 
\, + \, A(X \, , \, \Gamma)\, .
\edm
The vector $A(X,\Gamma)$ is defined by
\bdm
A(X \, , \, \Gamma) \ := \  \left\{\begin{array}{ll} 0& \mathrm{if} 
\, X \, \mathrm{and} \, \Gamma \, \mathrm{are \, proportional}\\ 
(n - 2) ||\Gamma||^2 \cdot X & \mathrm{if} \,  X \, \mathrm{and} \,
\Gamma \, \mathrm{are \,
 orthogonal}  \end{array} \right. 
\edm
The scalar curvature $\Scal^g$ is given by the formula
\bdm
\Scal^g \ = \ 2 \, (1 \, - \, n) \, \delta^g(\Gamma) \, + \,
(n \, - \, 1)(n \, - \, 2) \, ||\Gamma||^2  .
\edm
\end{thm}
\begin{proof}
The spinor field $\Psi$ is parallel with respect to the
connection $\nabla^{\mathrm{vec}}$, i.\,e. 
\bdm
\nabla^g_X \Psi \ = \ \frac{1}{2} \cdot
(X \wedge \Gamma) \cdot \Psi \, .
\edm
We compute the square of the Dirac operator as well as the
spinorial Laplacian on $\Psi$,
\begin{eqnarray*}
\big(D^g)^2 \Psi &=& \frac{1 - n}{2} \, \big( \delta^g(\Gamma) \, + \, 
d \Gamma \big) \cdot \Psi \, + \, \frac{(n-1)^2}{4} \, ||\Gamma||^2 \cdot 
\Psi \, , \\
\Delta\, \Psi &=& - \, \frac{1}{2} \cdot d \Gamma \cdot \Psi \, + \, 
\frac{n - 1}{2} \, ||\Gamma||^2 \cdot \Psi \, .
\end{eqnarray*}
The Schr\"odinger-Lichnerowicz formula $(D^g)^2 = \Delta + \Scal^g/4$
yields the equation
\bdm
2 \, (1 - n) \, \delta^g(\Gamma) \cdot \Psi \, + \, (n - 1)(n - 2) \, 
||\Gamma||^2 \cdot \Psi \, + \, 2 \, (2 - n) \, d \Gamma \cdot \Psi \ = \
\Scal^g \cdot \Psi \, .
\edm
Then we conclude that $d \Gamma \cdot \Psi = 0$ and
\bdm
2 \, (1 - n) \, \delta^g(\Gamma)  \, + \, (n - 1)(n - 2) \, 
||\Gamma||^2  \ = \ \Scal^g  .
\edm
The differential equation for the spinor $\Psi$
allows us to compute the action of the curvature
$\mathcal{R}^g(X,Y) \cdot \Psi =  \nabla^g_X \nabla^g_Y \Psi  -  
\nabla^g_Y \nabla^g_X \Psi  -  \nabla^g_{[X,Y]} \Psi$ on the spinor.
Then we use the well-known formula (see \cite{Fri1})
\bdm
\Ric^g(X) \cdot \Psi \ = \ - \, 2 \, \sum_{i=1}^n e_i \cdot \mathcal{R}^g(X ,
e_i) \cdot \Psi \ ,
\edm
and after a straightforward algebraic calculation
in the Clifford algebra we obtain
\bdm
\Ric^g(X) \cdot \Psi \ = \ A(X \, , \, \Gamma) \cdot \Psi \, + \,
(n - 3) \, (\nabla_X^g\Gamma) \cdot \Psi \, - \, \delta^g(\Gamma) \, X \cdot
\Psi \, 
+ \, \sum_{i=1}^n g(X \, , \, 
\nabla_{e_i}^g \Gamma) \cdot e_i \cdot \Psi \ .
\edm
Consider the inner product of the latter equation by the spinor
$Y \cdot \Psi$. Then we obtain
\bdm
g(\Ric^g(X),  Y) \, = \, g(A(X ,  \Gamma), Y)  +  
(n - 3) \, g(\nabla_X^g\Gamma , Y)  -  \delta^g(\Gamma) \, 
g(X ,  Y)  + g(X , \nabla^g_Y \Gamma) \, . 
\edm
Since the Riemannian Ricci tensor $\Ric^g$ is symmetric, the antisymmetric
part
\bdm
\frac{n - 4}{2} \, d \Gamma(X \, , \, Y)
\edm
of the right side has to vanish. The formula for the Ricci tensor
follows immediately.
\end{proof}
\begin{NB}
The conditions of the latter theorem are satisfied
for the groups $\G_2 \subset \SO(7)$ and $\Spin(7) \subset \SO(8)$.
The subgroups $\U(n) \subset \SO(2n)$ or 
$\Spin(9) \subset \SO(16)$ do {\it not} satisfy the
conditions, there are no invariant spinors.
\end{NB}
\vspace{2mm}

\begin{NB}
In dimension $n = 4$,
the condition $d \Gamma \cdot \Psi = 0$ defines
a $3$-dimensional subspace $V^3(\Psi) \subset \Lambda^2(\R^4)$ of $2$-forms
depending on the spinor $\Psi$. It is the isotropy Lie algebra
of the spinor $\Psi$.
\end{NB}
\vspace{2mm}

\noindent
For the special vector $X = \Gamma$, the formula for the Ricci tensor 
simplifies,
\bdm
\Ric^g(\Gamma)  \ = \ (n - 2) \, \nabla^g_{\Gamma}\Gamma 
\, - \,  \delta^g(\Gamma) \cdot \Gamma \, . 
\edm
We multiply the latter equation by the vector field $\Gamma$.  
In this way we  obtain the product $g(\Ric^g(\Gamma) \, , \, \Gamma)$
of the two vectors.
\vspace{2mm}

\begin{cor}
Suppose that the subgroup $\G \subset \SO(n)$ lifts into the spin group and
admits an invariant spinor $0 \neq \Psi \in \Delta_n$. 
Then, for any $\G$-structure of vectorial type, we have
\bdm
g( \mathrm{Ric}^g(\Gamma) \, , \, \Gamma) \ = \ \frac{(n -2)}{2} \cdot
\Gamma(||\Gamma||^2) \, - \, \delta^g(\Gamma) \cdot ||\Gamma||^2  .
\edm
If the manifold $M^n$ is compact, then
\bdm
\int_{M^n}g( \mathrm{Ric}^g(\Gamma) \, , \, \Gamma) \ = \ \frac{(n - 4)}{2} 
\cdot \int_{M^n} \delta^g(\Gamma) \cdot ||\Gamma||^2  .
\edm
\end{cor}
\vspace{2mm}

\noindent
The next Proposition states that -- up to a conformal change of the metric --
compact $\G$-structures of vectorial type admitting a parallel spinor
are locally conformal to products of Einstein spaces by $\R$. The
compactness is an essential assumption here. In \cite{ChiossiFino} the
authors constructed non-compact, $7$-dimensional
solvmanifolds equipped with a 
$\G_2$-structure of vectorial that are  not  Riemannian products of
$\R$ by an Einstein space of positive curvature.  
\vspace{2mm}

\begin{prop} \label{splitting}
Let $\G \subset \SO(n)$ be a subgroup that can be lifted into the spin 
group and
suppose that there exists a spinor $\G$-invariant $0 \neq \Psi \in \Delta_n$.
Consider a $\G$-structure of vectorial type
on a compact manifold and suppose that $\delta^g({\Gamma}) =  0$ holds. 
In dimension $n = 4$ we assume moreover that $\Gamma$ is a closed
form, $d \Gamma = 0$. Then we have
\begin{enumerate}
\item $\nabla^g \Gamma = 0$ .
\item $\mathrm{Ric}^g(\Gamma) = 0$.
\item If $X$ is orthogonal to $\Gamma$, then $\mathrm{Ric}^g(X) = 
(n - 1) \cdot ||\Gamma||^2 \cdot X$ .
\item The scalar curvature is positive
\bdm
\Scal^g \ = \ (n -  1)(n - 2) ||\Gamma||^2 \ > 0 \, .
\edm
\item The Lie derivative of any $\G$-invariant differential form
$\Sigma^k \in \Lambda^k(\R^n)$ vanishes
\bdm
\mathcal{L}_{\Gamma} \Sigma^k \ = = \nabla^g_{\Gamma} \Sigma^k \ = \ 0 \, .
\edm
\item The universal covering $\tilde{M}^n = Y^{n-1} \times \R^1$ splits
into $\R$ and an Einstein manifold  $Y^{n-1}$
with positive scalar curvature admitting a real Riemannian Killing spinor.
\end{enumerate}
\end{prop}
\begin{proof}
The $1$-form $\Gamma$ is by assumption harmonic and 
the Bochner formula for $1$-forms yields 
\bdm
0 \ = \ \int_{M^n} ||\nabla^{g}\Gamma||^2 \ + \ \frac{1}{3} 
\int_{M^n}g( \mathrm{Ric}^g(\Gamma) \, , \, \Gamma) \ =
\int_{M^n} ||\nabla^{g}\Gamma||^2 \, .
\edm
Consequently, $\Gamma$ is parallel with respect to the Levi-Civita
connection. Moreover, the restriction of the spinor field $\Psi$ to
the submanifold $Y^{n-1}$ defines a spinor field such that
\bdm
\nabla^{Y^{n-1}}_X \Psi \ = \ \frac{1}{2} \cdot X \cdot \Gamma \cdot \Psi \ ,
  \quad
\nabla^{Y^{n-1}}_X \Gamma \cdot \Psi \ = \ \frac{1}{2} ||\Gamma||^2 
\cdot X \cdot  \Psi 
\edm
holds for any vector $X \in T(Y^{n-1})$.
The spinor field $\Psi^* := \||\Gamma|| \cdot \Psi + \Gamma \cdot \Psi$
is a Killing spinor on $Y^{n-1}$.
\end{proof}
\begin{NB}
Let us discuss the latter proposition from the point of view
of Weyl geometry. 
Theorem \ref{Fixspinor} means that any $\G$-structure with a fixed 
spinor on a compact manifold
induces a Weyl-Einstein geometry with a closed
form $\Gamma$ ($n \geq 5$). Indeed, after a conformal change of the metric
the condition $\delta^{g}(\Gamma) = 0$ is satisfied. 
In this sense, Proposition \ref{splitting}
is a reformulation of Theorem 3 in \cite{Gauduchon}. In dimensions $n=7$ and
$n = 8$ this splitting has been discussed in \cite{Ivanov&Co}.
\end{NB}
%
\section{Geometric structures of vectorial type admitting a characteristic connection}\noindent
%
\noindent
Consider a geometric structure $\mathcal{R} \subset \mathcal{F}(M^{n})$
of vectorial type and suppose that it admits a characteristic connection.
Then the intrinsic torsion $\Gamma \in \R^n \otimes \mathg{m}$ is given by
a vector $\Gamma \in \R^n$,
\bdm
\Gamma (X) \ = \ \mathrm{pr}_{\m}(X \wedge \Gamma) \, .
\edm 
On the other side, there exists a $3$-form $\T^{\mathrm{c}}$ such that
\bdm
2 \cdot \Gamma(X) \ = \ - \, \Theta(\T^{\mathrm{c}})(X) \ 
= \ - \, \mathrm{pr}_{\m}(X \haken\T^{\mathrm{c}}) 
\edm
holds. Consequently, the vector field $\Gamma$ and the characteristic
torsion $\T^{\mathrm{c}}$ are related by the condition
\bdm
2 \cdot (X \wedge \Gamma) \, + \, X \haken \T^{\mathrm{c}} \ \in \, \g 
\edm
for all vectors $X$. In this case, we have two connection
$\nabla^{\mathrm{vec}}$
and $\nabla^{\mathrm{c}}$ preserving the $\G$-structure. The map
$\R^n \,\subset \R^n \otimes \m$ is injective for any subgroup
$\G \ne \SO(n)$. If, moreover, the map $\Theta : \Lambda^3(\R^n)
\rightarrow \R^n \otimes \m$ is injective too, then the characteristic
torsion $\T^{\mathrm{c}}$ is uniquely defined by the vector field
$\Gamma$. Structures with this property and with a non trivial $\Gamma$ 
cannot occur for all geometric structures. Indeed,
the $\G$-representation $\R^n$ has to be contained in the
$\G$-representation $\Lambda^3(\R^n)$. For example, for the subgroups
$\G = \SO(3) \subset \SO(5)$ , $\Spin(9) \subset 
\SO(16)$ or $\G = \mathrm{F}_4 \subset \SO(26)$
this condition is not satisfied (see \cite{Fri2}). In dimensions $n=7,8$ 
any $\G_2$- or $\Spin(7)$-structure of vectorial type 
admits a characteristic connection (see \cite{FriedrichIvanov}, \cite{Fri2}).
\vspace{2mm}

\noindent
If the group $\G$ preserves a spinor,  the
corresponding spinor
field $\Psi$ on the manifold is parallel with respect to the connections
$\nabla^{\mathrm{vec}}$ and $\nabla^{\mathrm{c}}$. A similar computation as 
in the proof of Theorem \ref{Fixspinor} yields the following formulas
linking $\Gamma$, $\T^{\mathrm{c}}$ and $\Psi$.
\vspace{2mm}

\begin{thm}
Let $\G \subset \SO(n)$ be a subgroup lifting into the spin group and
suppose that there exists a $\G$-invariant spinor  $0 \neq \Psi \in \Delta_n$.
Consider  a $\G$-structure of vectorial
type that admitas a characteristic connection. Denote by $\Gamma$ the
corresponding vector field and by $\T^{\mathrm{c}}$ the torsion of the
characteristic connection. Then we have
\begin{eqnarray*}
\big(\Gamma \haken \T^{\mathrm{c}}\big) \cdot \Psi &=&0 \, , 
\quad
\delta(\T^{\mathrm{c}}) \cdot \Psi \ = \ 0 \, , \quad
\T^{\mathrm{c}} \cdot \Psi \ = \ \frac{2}{3} \, (n - 1) \, \Gamma \cdot
\Psi \, , \\
(\T^{\mathrm{c}})^2 \cdot \Psi &=& \frac{4}{9} \, (n - 1)^2 
\, ||\Gamma||^2 \cdot \Psi \, , \\
d \T^{\mathrm{c}} \cdot \Psi &=& \frac{1}{3} \big( || \T^{\mathrm{c}}||^2 
\, - \, \frac{4}{9} \, (n - 1)^2 \, ||\Gamma||^2  \, - \, 
\Scal^{\nabla^{\T^{\mathrm{c}}}}
\big) \cdot \Psi \, , \\
2 \, (n - 1) \, \delta^g(\Gamma) &=& 2 \, \big( \, \frac{4}{9} \, (n - 1)^2 \, 
||\Gamma||^2 \, - \, ||\T^{\mathrm{c}}||^2 \big) \, - \, 
\Scal^{\nabla^{\T^{\mathrm{c}}}} \, .
\end{eqnarray*}
\end{thm}
\vspace{2mm}

\begin{exa} 
Consider a $7$-dimensional Riemannian manifold $(M^7 , g)$
equipped with a $\G_2$-structure of vectorial type, i.\,e., with
a generic $3$-form $\omega$. The differential equations
defining the vectorial type of the the structure read as
\bdm
d \, \omega \ = \ - \, 3 \, (\Gamma \wedge \omega) \, , \quad
\delta (\omega) \ = \ 4 \, ( \Gamma \haken \omega) \, . 
\edm
The characteristic torsion is given by the formula $\T^{\mathrm{c}} = 
- \, * (\Gamma \wedge \omega)$, see \cite{FriedrichIvanov}. In particular,
we have
\bdm
\Gamma \haken \T^{\mathrm{c}} \ = \ 0 \, , \quad
\delta(\T^{\mathrm{c}}) \ = \ 0 \, , \quad
||\T^{\mathrm{c}}||^2 \ = \ 4 \, ||\Gamma||^2 \, , \quad
12 \, \delta(\Gamma) \ = \ 6 \, ||\T^{\mathrm{c}}||^2 \, - \, 
\Scal^{\nabla^{\T^{\mathrm{c}}}} \, .
\edm
\end{exa}
%
%
\section{Generalized Hopf structures}\noindent

\noindent
The condition $\nabla^{\mathrm{vec}} \Gamma = 0$ or 
$\nabla^{\mathrm{vec}}\T^{\mathrm{c}} = 0$
is very restrictive. Indeed, it implies that
\bdm
\delta^g(\Gamma) \ = \ (n \, - \, 1) \cdot ||\Gamma||^2 .
\edm
Integrating the latter equation over a compact manifold, we obtain 
$\Gamma \equiv 0$. The conditions $\nabla^{\mathrm{c}} \Gamma = 0$ or
$\nabla^{\mathrm{c}}\T^{\mathrm{c}} = 0$ are more interesting
(see \cite{AlFrSchoe}).

\vspace{2mm}
\begin{prop} \label{Killing}
Suppose that $\Theta : \Lambda^3(\R^n)
\rightarrow \R^n \otimes \m$ is injective and let $\mathcal{R}$ be
a $\G$-structure of vectorial type admitting a characteristic
connection. If $\nabla^{\mathrm{c}}\Gamma = 0$,  then
\bdm
\delta^g(\Gamma) \ = \ 0 \, , \quad \delta^g(\T^{\mathrm{c}}) \ = \ 0 \, , 
\quad d \Gamma \ = \ \Gamma \haken \T^{\mathrm{c}} \, ,
\quad 2 \cdot \nabla^g \Gamma \ = \ d \Gamma \ .
\edm
In particular, $\Gamma$ is a Killing vector field.
\end{prop}
\begin{proof}
The formulas follow directly from the assumption,
\bdm
0\ = \ \nabla^{\mathrm{c}}_X \Gamma \ = \ \nabla^g_X \Gamma \, + \, \frac{1}{2}
\, \T^{\mathrm{c}}(X\, ,\, \Gamma \, , \, - ) \ . \qedhere
\edm
\end{proof}
\noindent
In complex geometry, a hermitian manifold of vectorial type such that
its characteristic torsion $\T^{\mathrm{c}}$ is
$\nabla^{\mathrm{c}}$-parallel
is called a {\it generalized Hopf manifold}.
These  $\mathcal{W}_4$-manifolds
have been studied by Vaisman, see \cite{Vaisman}. Let us revisit
this geometry in more detail.
\vspace{2mm}

\begin{exa}
Consider the subgroup $\U(2) \subset \SO(4)$. There are only two types
of $\U(2)$-structures. Moreover, a $\U(2)$-structure is of vectorial type
if and only if it admits a characteristic connection. The link between the
$3$-form $\T^{\mathrm{c}}$ and the vector field $\Gamma$ is $ \Gamma = 
* \, \T^{\mathrm{c}}$ (see \cite{AlFrSchoe}). Consequently, we obtain
\bdm
\nabla^{\mathrm{c}}_X \Gamma \ = \ \nabla^g_X \Gamma \, + \, 
\frac{1}{2} \, \T^{\mathrm{c}}(X \, , \, \Gamma \, , \, -) \ = \ 
\nabla^g_X \Gamma \, .
\edm
The condition $\nabla^{\mathrm{c}}_X \Gamma = 0$ is equivalent
to $\nabla^g_X \Gamma = 0$. These are generalized Hopf surfaces. 
They are locally conformal K\"ahler manifolds ($d \Gamma = 0$)
with a non-parallel vector field ($\nabla^g \Gamma \neq 0$ , see 
\cite{Belgun}).
There 
are also $\U(2)$-structures of vectorial type with a non-closed
form $\Gamma$ (see Example \ref{ExaU(2)}).
\end{exa}
\vspace{2mm}

\begin{exa}
Consider a hermitian manifold $(M^{6} ,  g ,  \J))$ and denote
by $\Omega$ its K\"ahler form. The vector field $\Gamma$ (the vector
part of the intrinsic torsion) is defined by
\bdm
\delta^g (\Omega) \ = \ 4 \cdot \J(\Gamma) \ = \ 4 \cdot (\Gamma \haken 
\Omega)
\, , \quad d \Omega \ = \ - \, 2 \cdot (\Gamma \wedge \Omega) \, .
\edm
Suppose that $M^{6}$ is of vectorial type and that 
$\nabla^{\mathrm{c}} \Gamma = 0$ holds.
Then its characteristic connection as well as the differential are given
by the formulas (see \cite{AlFrSchoe})
\bdm
\T^{\mathrm{c}} \ = \ 2 \cdot \big(\J(\Gamma) \wedge \Omega \big) \, ,  
\quad d \Gamma \ = \ 0 \ = \ \Gamma \haken \T^{\mathrm{c}} \, , 
\quad \nabla^g \Gamma \ = \ 0 \, .
\edm
In particular, $\Gamma$ is parallel with respect to the Levi-Civita 
connection and $\J(\Gamma)$ is a Killing vector field.
\end{exa} 
\vspace{2mm}

\begin{dfn}
A $\G$-structure $\mathcal{R} \subset \mathcal{F}(M^n)$ of vectorial type
and admitting a characteristic connection is called a 
{\it generalized Hopf structure} if
$\nabla^{\mathrm{c}} \Gamma = 0$ holds.
\end{dfn}
\vspace{2mm}

\noindent
The vector field $\Gamma$ of a Hopf $\G$-structure is a Killing vector
field. $\Gamma$ is  
$\nabla^g$-parallel if and only if $d \Gamma = 0$ holds. Proposition
\ref{dGamma} and Theorem \ref{Fixspinor} contain sufficient 
conditions that the vector field of any
Hopf $\G$-structure is $\nabla^g$-parallel. This situation
occurs for the standard geometries of the groups $\G = \G_2 , \Spin(7)$ 
and for $\U(n)\, , \, n \geq 3$. However, there are
subgroups $\G \subset \SO(n)$ and Hopf $\G$-structures 
($\nabla^{\mathrm{c}} \Gamma = 0$) with a non $\nabla^g$-parallel
vector field.
%
    
\end{document}